\documentclass[amstex,12pt,russian,amssymb]{article}

\usepackage{mathtext}
\usepackage[cp1251]{inputenc}
\usepackage[T2A]{fontenc}
\usepackage[russian]{babel}
\usepackage[dvips]{graphicx}
\usepackage{amsmath}
\usepackage{amssymb}
\usepackage{amsxtra}
\usepackage{latexsym}
\usepackage{ifthen}

\begin{document}

\newcounter{lemma}
\newcommand{\lemma}{\par \refstepcounter{lemma}%
{\bf Лема \arabic{lemma}.}}

\newcounter{corollary}
\newcommand{\corollary}{\par \refstepcounter{corollary}%
{\bf Наслідок \arabic{corollary}.}}

\newcounter{remark}
\newcommand{\remark}{\par \refstepcounter{remark}%
{\bf Зауваження \arabic{remark}.}}

\newcounter{theorem}
\newcommand{\theorem}{\par \refstepcounter{theorem}%
{\bf Теорема \arabic{theorem}.}}

\newcounter{proposition}
\newcommand{\proposition}{\par \refstepcounter{proposition}%
{\bf Твердження \arabic{proposition}.}}

\newcounter{example}
\newcommand{\example}{\par \refstepcounter{example}%
{\bf Приклад \arabic{example}.}}

\renewcommand{\refname}{\centerline{\bf Список літератури}}

\renewcommand{\figurename}{Мал.}

\newcommand{\proof}{{\it Доведення.\,\,}}

\noindent УДК 517.5

\medskip\medskip
{\bf Є.О.~Севостьянов} (Житомирський державний університет імені
Івана Фран\-ка; Інститут прикладної математики і механіки НАН
України, м.~Слов'янськ)

{\bf В.А.~Таргонський} (Житомирський державний університет імені
Івана Фран\-ка)

\medskip\medskip
{\bf E.O.~Sevost'yanov} (Zhytomyr Ivan Franko State University;
Institute of Applied Ma\-the\-ma\-tics and Mechanics of NAS of
Ukraine, Slov'yans'k)

{\bf V.A.~Targonskii} (Zhytomyr Ivan Franko State University)

\medskip
{\bf Про збіжність послідовності відображень з оберненою модульною
умовою до дискретного відображення}

{\bf On convergence of a sequence of mappings with inverse modulus
inequality to a discrete mapping}

\medskip\medskip
Досліджені відображення, які задовольняють обернену нерівність типу
Полецького в області евклідового простору. Доведено, що рівномірна
границя сім'ї таких відображень є дискретним відображенням. Окремо
розглянуті області, локально зв'язні на своїй межі, та регулярні
області в квазіконформному сенсі.

\medskip\medskip
We have studied the mappings that satisfy the Poletsky-type inverse
inequality in the domain of the Euclidean space. It is proved that
the uniform boundary of the family of such mappings is a discrete
mapping. We separately considered domains that are locally connected
at their boundary and regular domains in the quasiconformal sense.

\newpage
{\bf 1. Вступ.} Дана стаття присвячена дослідженню відображень з
обмеженим і скінченним спотворенням, див., напр., \cite{Cr},
\cite{MRV$_1$}, \cite{MRSY}, \cite{Vu} і~\cite{Va}. У нашій
нещодавній публікації \cite{Sev$_2$} ми довели дискретність
відображення в замиканні області, яке задовольняє одну (обернену)
умову спотворення модуля сімей кривих та само є відкритим,
дискретним і замкненим у ній. Отримання цього результату передбачає
можливість неперервного продовження на межу, що за вказаних умов
також виконується (див. там же; див. також~\cite{SSD}). Дана замітка
присвячена отриманню більш загального результату: ми покажемо, що
аналогічну властивість задовольняють не тільки самі відображення,
але й рівномірні границі послідовності відображень. Зауважимо, що
результати, близькі до наведених у цій роботі, отримані для
квазірегулярних відображень в~\cite{Vu}.

\medskip
Звернемося до означень. Всюди далі $M(\Gamma)$ позначає модуль сім'ї
$\Gamma$ (див.~\cite{Va}). Нехай $y_0\in {\Bbb R}^n,$
$0<r_1<r_2<\infty$ і
\begin{equation}\label{eq1**}
A=A(y_0, r_1,r_2)=\left\{ y\,\in\,{\Bbb R}^n:
r_1<|y-y_0|<r_2\right\}\,.\end{equation}
Нехай $x_0\in{\Bbb R}^n,$ тоді покладемо
$$B(x_0, r)=\{x\in {\Bbb R}^n: |x-x_0|<r\}\,, \quad {\Bbb B}^n=B(0, 1)\,,$$
$$S(x_0,r) = \{
x\,\in\,{\Bbb R}^n : |x-x_0|=r\}\,. $$
Для заданих множин $E,$ $F\subset\overline{{\Bbb R}^n}$ і області
$D\subset {\Bbb R}^n$ позначимо через $\Gamma(E,F,D)$ сім'ю всіх
кривих $\gamma:[a,b]\rightarrow \overline{{\Bbb R}^n}$ таких, що
$\gamma(a)\in E,\gamma(b)\in\,F$ і $\gamma(t)\in D$ при $t \in [a,
b].$ Якщо $f:D\rightarrow {\Bbb R}^n$ -- задане відображення,
$y_0\in f(D)$ і $0<r_1<r_2<d_0=\sup\limits_{y\in f(D)}|y-y_0|,$ то
через $\Gamma_f(y_0, r_1, r_2)$ ми позначимо сім'ю всіх кривих
$\gamma$ в області $D$ таких, що $f(\gamma)\in \Gamma(S(y_0, r_1),
S(y_0, r_2), A(y_0,r_1,r_2)).$ Нехай $Q:{\Bbb R}^n\rightarrow [0,
\infty]$ -- вимірна за Лебегом функція, $Q(x)\equiv 0$ при $x\in
{\Bbb R}^n\setminus f(D).$ Будемо говорити, що {\it $f$ задовольняє
обернену нерівність Полецького} в точці $y_0\in
\overline{f(D)}\setminus\{\infty\},$ якщо співвідношення
\begin{equation}\label{eq2*A}
M(\Gamma_f(y_0, r_1, r_2))\leqslant \int\limits_{A(y_0,r_1,r_2)\cap
f(D)} Q(y)\cdot \eta^n (|y-y_0|)\, dm(y)
\end{equation}
виконується для довільної вимірної за Лебегом функції $\eta:
(r_1,r_2)\rightarrow [0,\infty ]$ такій, що
\begin{equation}\label{eqA2}
\int\limits_{r_1}^{r_2}\eta(r)\, dr\geqslant 1\,.
\end{equation}
За допомогою інверсії $\psi(y)=\frac{y}{|y|^2}$ можна також
визначити співвідношення~(\ref{eq2*A}) у точці $y_0=\infty.$
Відображення $f:D\rightarrow {\Bbb R}^n$ називається {\it
дискретним}, якщо прообраз $\{f^{-1}\left(y\right)\}$ кожної точки
$y\,\in\,{\Bbb R}^n$ складається з ізольованих точок, і {\it
відкритим}, якщо образ будь-якої відкритої множини $U\subset D$ є
відкритою множиною в ${\Bbb R}^n.$ Відображення $f$ області $D$ на
$D^{\,\prime}$ називається {\it замкненим}, якщо $f(E)$ є замкненим
в $D^{\,\prime}$ для будь-якої замкненої множини $E\subset D$ (див.,
напр., \cite[розд.~3]{Vu}). У подальшому,  в розширеному просторі
$\overline{{{\Bbb R}}^n}={{\Bbb R}}^n\cup\{\infty\}$
використовується {\it сферична (хордальна) метрика}
$h(x,y)=|\pi(x)-\pi(y)|,$ де $\pi$ -- стереографічна проекція
$\overline{{{\Bbb R}}^n}$  на сферу
$S^n(\frac{1}{2}e_{n+1},\frac{1}{2})$ в ${{\Bbb R}}^{n+1},$ а саме,
$$h(x,\infty)=\frac{1}{\sqrt{1+{|x|}^2}}\,,$$
\begin{equation}\label{eq3C}
\ \ h(x,y)=\frac{|x-y|}{\sqrt{1+{|x|}^2} \sqrt{1+{|y|}^2}}\,, \ \
x\ne \infty\ne y
\end{equation}
(див., напр., \cite[означення~12.1]{Va}). Всюди далі межа $\partial
A$ множини $A$ і замикання $\overline{A}$ слід розуміти в сенсі
розширеного евклідового простору $\overline{{\Bbb R}^n}.$ Неперервне
продовження відображення $f:D\rightarrow {\Bbb R}^n$ також (якщо
непорозуміння є неможливим) слід розуміти з точки зору відображення
зі значеннями в $\overline{{\Bbb R}^n}$ і відносно метрики $h$
у~(\ref{eq3C}). Нагадаємо, що область $D\subset {\Bbb R}^n$
називається {\it локально зв'язною в точ\-ці} $x_0\in\partial D,$
якщо для будь-якого околу $U$ точки $x_0$ знайдеться окіл $V\subset
U$ точки $x_0$ такий, що $V\cap D$ є зв'язним. Область $D$ локально
зв'язна на $\partial D,$ якщо $D$ локально зв'язна в кожній точці
$x_0\in\partial D.$ Ми кажемо, що область $D$ є {\it скінченно
зв'язною} у точці $x_0\in\partial D$, якщо для будь-якого околу $U$
точки $x_0$ існує окіл $V\subset U$ цієї точки такий, що $V\cap D$
складається зі скніченної кількості компонент зв'язності.
Відображення $f:X\rightarrow Y$ називається {\it нульвимірним,} якщо
для кожної точки $y\in Y$ множина $f^{\,-1}(y)$ не містить в собі
континуум $K\subset X,$ відмінний від точки. Межа області $D$
називається {\it слабко плоскою} в точці $x_0\in
\partial D,$ якщо для кожного $P>0$ і для будь-якого околу $U$
точки $x_0$ знайдеться окіл $V\subset U$ цієї ж самої точки такий,
що $M(\Gamma(E, F, D))>P$ для будь-яких континуумів $E, F\subset D,$
які перетинають $\partial U$ і $\partial V.$ Межа області $D$
називається слабко плоскою, якщо відповідна властивість виконується
в будь-якій точці межі $D.$

\medskip
Для ві\-доб\-ра\-жен\-ня $f:D\,\rightarrow\,{\Bbb R}^n,$ множини
$E\subset D$ та $y\,\in\,{\Bbb R}^n,$ визначимо {\it функцію
кратності $N(y,f,E)$} як кількість прообразів точки $y$ у множині
$E,$ тобто
$$
N(y,f,E)\,=\,{\rm card}\,\left\{x\in E: f(x)=y\right\}\,,
$$
$$
N(f,E)\,=\,\sup\limits_{y\in{\Bbb R}^n}\,N(y,f,E).
$$
Будемо говорити, що функція ${\varphi}:{\Bbb R}^n\rightarrow{\Bbb
R}$ має {\it скінченне середнє коливання} у точці $x_0\in D,$ пишемо
$\varphi\in FMO(x_0),$ якщо
$
\limsup\limits_{\varepsilon\rightarrow
0}\frac{1}{\Omega_n\varepsilon^n}\int\limits_{B( x_0,\,\varepsilon)}
|{\varphi}(x)-\overline{{\varphi}}_{\varepsilon}|\ dm(x)<\infty\,,
$
де $\overline{{\varphi}}_{\varepsilon}=\frac{1}
{\Omega_n\varepsilon^n}\int\limits_{B(x_0,\,\varepsilon)}
{\varphi}(x) dm(x).$ Покладемо
\begin{equation}\label{eq12}
q_{y_0}(r)=\frac{1}{\omega_{n-1}r^{n-1}}\int\limits_{S(y_0,
r)}Q(y)\,d\mathcal{H}^{n-1}(y)\,, \end{equation} де $\omega_{n-1}$
позначає площу одиничної сфери ${\Bbb S}^{n-1}$ в ${\Bbb R}^n,$
Виконуються наступні твердження.

\medskip
\begin{theorem}\label{th1}
{\sl Нехай $D\subset {\Bbb R}^n,$ $n\geqslant 2,$ -- область, яка
має слабо плоску межу, а область $D^{\,\prime}\subset {\Bbb R}^n$ є
скінченно зв'язною на своїй межі. Припустимо, $f_m,$ $m=1,2,\ldots
,$ -- послідовність відкритих, дискретних і замкнених відображень
області $D$ на $D^{\,\prime},$ які збігаються рівномірно у $D$ до
деякого (не тотожно сталого) відображення $f:D\rightarrow {\Bbb
R}^n$ і задовольняють співвідношення~(\ref{eq2*A}) в кожній точці
$y_0\in \overline{D^{\,\prime}}.$ Нехай також виконана одна з двох
умов:

\medskip
1) $Q\in FMO(\overline{D^{\,\prime}});$

\medskip
2) для будь-якого $y_0\in \overline{D^{\,\prime}}$ знайдеться
$\delta(y_0)>0$ таке, що для достатньо малих $\varepsilon>0$
виконуються умови
\begin{equation}\label{eq5B}
\int\limits_{\varepsilon}^{\delta(y_0)}
\frac{dt}{tq_{y_0}^{\frac{1}{n-1}}(t)}<\infty, \qquad
\int\limits_{0}^{\delta(y_0)}
\frac{dt}{tq_{y_0}^{\frac{1}{n-1}}(t)}=\infty\,.
\end{equation}
Тоді відображення $f$ є нульвимірним і замкненим у $D,$ крім того,
$f$ має неперервне продовження $\overline{f}:\overline{D}\rightarrow
\overline{D^{\,\prime}},$
$\overline{f}(\overline{D})=\overline{D^{\,\prime}}$ і
$\overline{f}$ є нульвимірним відносно $\overline{D}.$}
\end{theorem}

\medskip
У випадку, коли відображена область є не скінченно зв'язною, а
навіть локально зв'язною на своїй межі, крім того, $f$ зберігає
орієнтацію, маємо більш посилений варіант теореми~\ref{th1}.

\medskip
\begin{theorem}\label{th3}
{\sl Нехай в умовах теореми~\ref{th1} $D^{\,\prime}$ є локально
зв'язною на своїй межі, а відображення $f$ зберігає орієнтацію. Тоді
відображення $f$ є відкритим, дискретним і замкненим у $D,$ і має
неперервне продовження $\overline{f}:\overline{D}\rightarrow
\overline{D^{\,\prime}}$ таке, що $f(D)=D^{\,\prime},$
$f(\overline{D})=\overline{D^{\,\prime}}$ і $N(f, D)=N(f,
\overline{D})<\infty.$ Зокрема, $\overline{f}$ є дискретним у
$\overline{D}.$}
\end{theorem}

\medskip
{\bf 2. Основні леми.} Наступне твердження доведено
в~\cite[лема~3.1]{Sev$_2$} для окремого відображення, див. також
відповідний класичний варіант для квазірегулярних відображень
у~\cite[лема~4.4]{Vu}. Версія, яка наведена нижче, стосується
збіжної послідовності відображень.

\medskip
\begin{lemma}\label{lem1}
{\sl\, Нехай $D\subset {\Bbb R}^n,$ $n\geqslant 2,$ -- область, яка
має слабо плоску межу. Припустимо, $f_m,$ $m=1,2,\ldots ,$ --
послідовність відкритих дискретних відображень області $D$ на
$D^{\,\prime},$ які збігаються рівномірно у $D$ до деякого (не
тотожно сталого) відображення $f:D\rightarrow {\Bbb R}^n$ і
задовольняють співвідношення~(\ref{eq2*A}) в кожній точці $y_0\in
\overline{D^{\,\prime}}.$ Припустимо, що для кожного $y_0\in
\overline{D^{\,\prime}}$ знайдеться
$\varepsilon_0=\varepsilon_0(y_0)>0$ і вимірна за Лебегом функція
$\psi:(0, \varepsilon_0)\rightarrow [0,\infty]$ такі, що
\begin{equation}\label{eq7***} I(\varepsilon,
\varepsilon_0):=\int\limits_{\varepsilon}^{\varepsilon_0}\psi(t)\,dt
< \infty\quad \forall\,\,\varepsilon\in (0, \varepsilon_0)\,,\quad
I(\varepsilon, \varepsilon_0)\rightarrow
\infty\quad\text{при}\quad\varepsilon\rightarrow 0\,,
\end{equation}
і, крім того, при $\varepsilon\rightarrow 0$
\begin{equation} \label{eq3.7.2}
\int\limits_{A(y_0, \varepsilon, \varepsilon_0)}
Q(y)\cdot\psi^{\,n}(|y-y_0|)\,dm(x) = o(I^n(\varepsilon,
\varepsilon_0))\,,\end{equation}
де $A(y_0, \varepsilon, \varepsilon_0)$ визначено в (\ref{eq1**}).
Нехай $C_j,$ $j=1,2,\ldots ,$ -- довільна послідовність континуумів
така, що $h(C_j)\geqslant \delta>0$ для деякого $\delta>0$ і всіх
$j\in {\Bbb N}$ і, крім того, $h(f(C_j))\rightarrow 0$ при
$j\rightarrow\infty.$ Тоді для кожного $y_0\in
\overline{D^{\,\prime}}$ знайдеться $\delta_1>0$ таке, що
\begin{equation}\label{eq1}
h(f(C_j), y_0)\geqslant \delta_1>0
\end{equation}
для всіх $m\in {\Bbb N}$ і $j\in {\Bbb N}.$}
\end{lemma}

\medskip
\begin{remark}
У точці $y_0=\infty$ співвідношення~(\ref{eq3.7.2}) слід розуміти за
допомогою інверсії $\psi(y)=\frac{y}{|y|^2}$ в нулі. Тобто, замість
$$\int\limits_{A(y_0, \varepsilon, \varepsilon_0)}
Q(y)\cdot\psi^{\,n}(|y-y_0|)\,dm(y) = o(I^n(\varepsilon,
\varepsilon_0))$$
буде розглядатися умова
$$\int\limits_{A(0, \varepsilon, \varepsilon_0)}
Q\left(\frac{y}{|y|^2}\right)\cdot\psi^{\,n}(|y|)\,dm(y) =
o(I^n(\varepsilon, \varepsilon_0))\,.$$
\end{remark}

\medskip
{\it Доведення леми~\ref{lem1}}. Надалі вважаємо $y_0\ne\infty.$
Припустимо супротивне, а саме, нехай $h(f(C_{j_k}), y_0)\rightarrow
0$ при $k\rightarrow\infty$ для деякої зростаючої послідовності
номерів $j_k,$ $k=1,2,\ldots .$ З огляду
на~\cite[лема~1]{SevSkv$_3$} $f$ або нульвимірне, або стале
відображення.  Оскільки за припущенням $f$ не є сталим, $f$ --
нульвимірне. Тоді існує точка $z_0\in D$ така, що $f(z_0)\ne y_0.$
Оскільки $f$ неперервне відображення, можна знайти $r_0>0$ таке, що
\begin{equation}\label{eq3}
\overline{B(z_0, r_0)}\subset D\,,\qquad f(\overline{B(z_0,
r_0)})\cap \overline{B(y_0, \varepsilon_1)}=\varnothing\,,
\end{equation}
де $\varepsilon_1>\varepsilon_0$ і $\varepsilon_0$ -- число зі
співвідношень~(\ref{eq7***})--(\ref{eq3.7.2}).
Покладемо $F:=\overline{B(z_0, r_0)}$ і позначимо
$\Gamma_k:=\Gamma(F, C_{j_k}, D).$

З огляду на~(\ref{eq3}) покажемо, що існує окіл $V$ елементу $y_0$
такий, що
\begin{equation}\label{eq5A}
f_j(F)\cap V=\varnothing\,,\qquad j=1,2,\ldots\,.
\end{equation}
Припустимо протилежне. Тоді існують послідовності $j_m\in{\Bbb N}$ і
$x_m\in F,$ $m=1,2,\ldots, $ такі, що $f_{j_m}(x_m)\rightarrow y_0$
при $m\rightarrow\infty.$ З огляду на те, що $f_{j_m}(F)$ є
компактом у $D^{\,\prime}$ при кожному фіксованому $m,$
послідовність індексів $j_m,$ $m=1,2,\ldots,$ можна вважати
зростаючою. Тоді за нерівністю трикутника
\begin{equation}\label{eq5F}
h(f(x_m), y_0)\leqslant h(f(x_m), f_{j_m}(x_m))+ h(f_{j_m}(x_m),
y_0)\rightarrow 0
\end{equation}
при $m\rightarrow\infty,$ оскільки $h(f(x_m),
f_{j_m}(x_m))\rightarrow 0$ при $m\rightarrow\infty$ по рівномірній
збіжності $f_m$ до $f$ у $D$ по метриці $h.$
Співвідношення~(\ref{eq5F}) суперечить~(\ref{eq3}), що
доводить~(\ref{eq5A}).

\medskip
Зауважимо, що область зі слабко плоскою межею є рівномірною (див.
\cite[наслідок~4.3]{Sev$_2$}), тобто, якщо континууми $C_{j_k}$
задовольняють умову $h(C_{j_k})\geqslant \delta>0$ і
$h(F)>\delta_*>0,$ то знайдеться $\delta_2>0$ таке, що
\begin{equation}\label{eq3A}
M(\Gamma_k)\geqslant \delta_2>0
\end{equation}
для всіх $k\in {\Bbb N}.$ Доведемо, що для кожного $l\in {\Bbb N}$
існує номер $k=k_l$ такий, що
\begin{equation}\label{eq3B}
f_{j_k}(C_{j_k})\subset B(y_0, 1/l)\,,\qquad k\geqslant k_l\,.
\end{equation}
Припустимо протилежне. Тоді існує $l_0\in {\Bbb N}$ таке, що
\begin{equation}\label{eq3F}
f_{j_{k_m}}(C_{j_{k_m}})\cap ({\Bbb R}^n\setminus B(y_0,
1/l_0))\ne\varnothing
\end{equation}
для деякої зростаючої послідовності номерів $k_m,$ $m=1,2,\ldots .$
В такому випадку, знайдеться послідовність $z_{m}\in
f_{j_{k_m}}(C_{j_{k_m}})\cap ({\Bbb R}^n\setminus B(y_0, 1/l_0)),$
$m\in {\Bbb N}.$ Нехай $z_m=f_{j_{k_m}}(x_m),$ $x_m\in C_{j_{k_m}},$
$m=1,2,\ldots .$ Оскільки за припущенням $h(f(C_{j_{k_m}}),
y_0)\rightarrow 0$ для деякої зростаючої послідовності номерів
$j_k,$ $k=1,2,\ldots ,$ то зокрема
\begin{equation}\label{eq3E}
h(f(C_{j_{k_m}}), y_0)\rightarrow 0\qquad {\text при}\qquad
m\rightarrow\infty\,.
\end{equation}
Оскільки $h(f(C_{j_{k_m}}), y_0)=\inf\limits_{y\in
f(C_{j_{k_m}})}h(y, y_0)$ і $f(C_{j_{k_m}})$ є компактом як
неперервний образ компакту $C_{j_{k_m}}$ при відображенні $f,$
звідси випливає, що $h(f(C_{j_{k_m}}), y_0)=h(y_m, y_0),$ де $y_m\in
f(C_{j_{k_m}}).$ Тоді зі співвідношення~(\ref{eq3E}) ми отримаємо,
що $y_m\rightarrow y_0$ при $m\rightarrow\infty.$ Оскільки за умовою
$h(f(C_j))=\sup\limits_{y,z\in f(C_j)}h(y,z)\rightarrow 0$ при
$j\rightarrow\infty,$ то $h(y_m, f(x_m))\leqslant
h(f(C_{j_{k_m}}))\rightarrow 0$ при $l\rightarrow\infty.$ Тоді за
нерівністю трикутника і за доведеним вище ми будемо мати, що:
$$h(z_m, y_0)=h(f_{j_{k_m}}(x_m), y_0)\leqslant$$
\begin{equation}\label{eq8}
\leqslant h(f_{j_{k_m}}(x_m), f(x_m))+h(f(x_m), y_m)+h(y_m,
y_0)\rightarrow 0\,,\qquad m\rightarrow\infty\,,
\end{equation}
оскільки $h(f_{j_{k_m}}(x_m), f(x_m))\leqslant\sup\limits_{x\in D}
h(f_{j_{k_m}}(x), f(x))\rightarrow 0$ при $m\rightarrow\infty$ за
рівномірною збіжністю $f_j$ до $f$ при $j\rightarrow\infty;$ останні
два члени в~(\ref{eq8}) також прямують до нуля за побудовою та
доведеним вище. Зі співвідношення~(\ref{eq8}) випливає, що
$z_m\rightarrow y_0$ при $m\rightarrow\infty,$ що суперечить
визначенню $z_m.$ Тоді припущення в~(\ref{eq3F}) є невірним, що
доводить~(\ref{eq3B}).

\medskip
З~(\ref{eq3B}) випливає, що
\begin{equation}\label{eq3M}
f_{j_{k_l}}(C_{j_{k_l}})\subset B(y_0, 1/l)\,,
\end{equation}
причому можна вважати послідовність номерів $k_l,$ $l=1,2,\ldots ,$
зростаючою. Переходячи у~(\ref{eq3M}) до перенумерації можна
вважати, що сама послідовність $j_k,$ $k=1,2, \ldots ,$ задовольняє
цю умову, тобто,
\begin{equation}\label{eq3N}
f_{j_k}(C_{j_{k}})\subset B(y_0, 1/k)\,, k=1,2,\ldots\,.
\end{equation}
Зменшуючи за потреби число $\varepsilon_0>0$ з умови леми можна
вважати, що окіл $V$ точки $y_0$ у~(\ref{eq5A}) є таким, що
$\overline{B(y_0, \varepsilon_0)}\subset V.$ Крім того, нехай
$k_0\in {\Bbb N}$ є таким, що $B(y_0, 1/k)\subset B(y_0,
\varepsilon_0)$ при всіх $k\geqslant k_0.$
В такому випадку, зауважимо, що
\begin{equation}\label{eq3G}
f_{j_{k}}(\Gamma_{j_{k}})>\Gamma(S(y_0, 1/k), S(y_0, \varepsilon_0),
A(y_0, 1/k,\varepsilon_0))\,.
\end{equation}
Дійсно, нехай $\widetilde{\gamma}\in f_{j_{k}}(\Gamma_{j_{k}}).$
Тоді $\widetilde{\gamma}(t)=f_{j_{k}}(\gamma(t)),$ де $\gamma\in
\Gamma_{j_{k}},$ $\gamma:[0, 1]\rightarrow D,$ $\gamma(0)\in F,$
$\gamma(1)\in C_{j_{k}}.$ За співвідношенням~(\ref{eq5A}) маємо, що
$f_{j_{k}}(\gamma(0))\in f(F)\subset {\Bbb R}^n\setminus
\overline{B(y_0, \varepsilon_0)},$ крім того, за
співвідношенням~(\ref{eq3N}) ми отримаємо, що
$f_{j_{k}}(\gamma(1))\in f_{j_{k}}(C_{j_{k}})\subset B(y_0,
\varepsilon_0).$ Отже, $|f(\gamma(t))|\cap B(y_0,
\varepsilon_0)\ne\varnothing \ne |f(\gamma(t))|\cap ({\Bbb
R}^n\setminus B(y_0, \varepsilon_0)).$ Тоді з огляду на
\cite[теорема~1.I.5.46]{Ku} ми отримаємо, що існує $0<t_1<1$ таке,
що $f_{j_{k}}(\gamma(t_1))\in S(y_0, \varepsilon_0).$ Покладемо
$\gamma_1:=\gamma|_{[t_1, 1]}.$ Можна вважати, що
$f_{j_{k}}(\gamma(t))\in B(y_0, \varepsilon_0)$ при всіх $t\geqslant
t_1.$ Міркуючи аналогічно, ми отримаємо точку $t_2\in (t_1, 1]$
таку, що $f_{j_{k}}(\gamma(t_2))\in S(y_0, 1/k).$ Покладемо
$\gamma_2:=\gamma|_{[t_1, t_2]}.$ Можна вважати, що
$f_{j_{k}}(\gamma(t))\not\in B(y_0, 1/k)$ при всіх $t\in [t_1,
t_2].$ Тоді крива $f_{j_{k}}(\gamma_2)$ є підкривою кривої
$f_{j_{k}}(\gamma)=\widetilde{\gamma},$ яка належить до сім'ї
$\Gamma(S(y_0, 1/k), S(y_0, \varepsilon_0), A(y_0,
1/k,\varepsilon_0)).$ Співвідношення~(\ref{eq3G}) встановлено.

\medskip
З~(\ref{eq3G}) випливає, що
\begin{equation}\label{eq3H}
\Gamma_{j_{k}}>\Gamma_{f_{j_{k}}}(S(y_0, 1/k), S(y_0,
\varepsilon_0), A(y_0, 1/k,\varepsilon_0))\,.
\end{equation}
Покладемо
$$\eta_{k}(t)=\left\{
\begin{array}{rr}
\psi(t)/I(1/k, \varepsilon_0), & t\in (1/k, \varepsilon_0)\,,\\
0,  &  t\not\in (1/k, \varepsilon_0)\,,
\end{array}
\right. $$
де $I(1/k, \varepsilon_0)=\int\limits_{1/k}^{\varepsilon_0}\,\psi
(t)\, dt.$ Зауважимо, що
$\int\limits_{1/k}^{\varepsilon_0}\eta_{l}(t)\,dt=1.$ Тоді за
співвідношеннями~(\ref{eq3.7.2}) і~(\ref{eq3H}), а також з огляду на
означення відображення $f_{j_{k}}$ у~(\ref{eq2*A}) будемо мати:
$$M(\Gamma_{j_{k}})\leqslant M(\Gamma_{f_{j_{k}}}(S(y_0, 1/k), S(y_0,
\varepsilon_0), A(y_0,1/k,\varepsilon_0)))\leqslant$$
\begin{equation}\label{eq3J}
\leqslant \frac{1}{I^n(1/k, \varepsilon_0)}\int\limits_{A(y_0, 1/k,
\varepsilon_0)} Q(y)\cdot\psi^{\,n}(|y-y_0|)\,dm(y)\rightarrow
0\quad \text{при}\quad m\rightarrow\infty\,.
\end{equation}
Останнє співвідношення суперечить~(\ref{eq3A}). Отримана
суперечність доводить лему.~$\Box$

\medskip
Ми кажемо, що область $D$ є {\it скінченно зв'язною} у точці
$x_0\in\partial D$, якщо для будь-якого околу $U$ точки $x_0$ існує
окіл $V\subset U$ цієї точки такий, що $V\cap D$ складається зі
скніченної кількості компонент зв'язності. Відображення
$f:X\rightarrow Y$ називається {\it нульвимірним,} якщо для кожної
точки $y\in Y$ множина $f^{\,-1}(y)$ не містить в собі континуум
$K\subset X,$ відмінний від точки. Наступна лема узагальнює
наслідок~4.5 у~\cite{Vu} на випадок відображень з необмеженою
характеристикою. Міркування, застосовані при її доведенні, по
більшій мірі аналогічні доведенню леми~3.2 в~\cite{Sev$_2$}.

\medskip
\begin{lemma}\label{lem2}
{\sl\, Нехай $D\subset {\Bbb R}^n,$ $n\geqslant 2,$ -- область, яка
має слабо плоску межу, а область $D^{\,\prime}\subset {\Bbb R}^n$ є
скінченно зв'язною на своїй межі. Припустимо, $f_m,$ $m=1,2,\ldots
,$ -- послідовність відкритих, дискретних і замкнених відображень
області $D$ на $D^{\,\prime},$ які збігаються рівномірно у $D$ до
деякого відображення $f:D\rightarrow {\Bbb R}^n$ і задовольняють
співвідношення~(\ref{eq2*A}) в кожній точці $y_0\in
\overline{D^{\,\prime}}.$ Нехай для будь-якого $y_0\in
\overline{D^{\,\prime}}$ знайдеться
$\varepsilon_0=\varepsilon_0(y_0)>0$ і вимірна за Лебегом функція
$\psi:(0, \varepsilon_0)\rightarrow [0,\infty]$ такі, що виконуються
співвідношення~(\ref{eq7***})--(\ref{eq3.7.2}), де $A(y_0,
\varepsilon, \varepsilon_0)$ визначено в (\ref{eq1**}). Тоді
відображення $f$ є нульвимірним і замкненим у $D,$ крім того, $f$
має неперервне продовження $\overline{f}:\overline{D}\rightarrow
\overline{D^{\,\prime}},$
$\overline{f}(\overline{D})=\overline{D^{\,\prime}}$ і
$\overline{f}$ є нульвимірним відносно $\overline{D}.$}
\end{lemma}

\medskip
\begin{proof}
Можна вважати $y_0\ne \infty.$ З огляду на \cite{SevSkv$_3$}
відображення $f$ є нульвимірним у $D.$ За лемою~5.1 в
\cite{Sev$_2$}, кожне $f_m,$ $m=1,2,\ldots ,$ має неперервне
продовження $\overline{f_m}:\overline{D}\rightarrow
\overline{D^{\,\prime}}$ на $\overline{D}$ таке, що
$f_m(\overline{D})=\overline{D^{\,\prime}}.$ Покажемо, що існує
$\lim\limits_{x\rightarrow x_0}f(x)$ для кожного $x_0\in\partial D.$
Дійсно, за нерівністю трикутника
$$h(f(x), f(x_0))\leqslant h(f(x), f_m(x))+h(f_m(x), f_m(x_0))+h(f_m(x_0), f(x_0))\,.$$
У цій нерівності перший і третій доданки прямують до нуля з огляду
на рівномірну збіжність $f_m$ у $D.$ Отже, для будь-якого
$\varepsilon>0$ знайдеться номер $M=M(\varepsilon)$ такий, що при
$m\geqslant M(\varepsilon)$
\begin{equation}\label{eq5D}
h(f(x), f(x_0))\leqslant \frac{2\varepsilon}{3}+h(f_m(x),
f_m(x_0))\,.
\end{equation}
При $m=M(\varepsilon)$ відображення $f_m$ неперервно продовжується в
точку $x_0,$ так що знайдеться $\delta=\delta(\varepsilon, x_0)>0$
таке, що $h(f_M(x), f_M(x_0))<\varepsilon/3$ при $0<|x-x_0|<\delta.$
Тоді з огляду на~(\ref{eq5D}) $h(f(x), f(x_0))<\varepsilon$ при тих
же $\delta.$

\medskip
Покажемо, що $f$ -- замкнене у $D.$ Нехай $E$ -- замкнене в $D.$
Треба показати, що $f(E)$ замкнене в $D^{\,\prime}.$ Нехай $z_m,$
$m=1,2,\ldots,$ -- послідовність у $f(E)$ така, що $z_m\rightarrow
z_0\in D^{\,\prime}$ при $m\rightarrow \infty.$ Покажемо, що $z_0\in
f(E).$ Оскільки $z_m\in f(E),$ знайдуться $x_m\in E$ такі, що
$f(x_m)=z_m,$ $m=1,2,\ldots .$ З огляду на компактність
$\overline{{\Bbb R}^n}$ можна вважати, що послідовність $x_m$
збігається до деякого $x_0\in \overline{D}$ при
$m\rightarrow\infty.$ Якщо $x_0\in D,$ то $x_0\in E,$ отже,
$f(x_0)=z_0\in E$ з огляду на неперервність відображення $f$ у $D.$
Якщо ж $x_0\in \partial D,$ то, оскільки $f_m$ має неперервне
продовження в точку $x_0,$ знайдеться послідовність $w_m\in D,$
$w_m\rightarrow x_0$ при $m\rightarrow\infty,$ така що $h(f_m(w_m),
f_m(x_0))<\frac{1}{m}.$ Тоді за нерівністю трикутника
$$h(f_m(x_0), f(x_0))\leqslant h(f_m(x_0), f_m(w_m))+h(f_m(w_m), f(w_m))+
h(f(w_m), f(x_0))\leqslant$$
\begin{equation}\label{eq6}\leqslant\frac{1}{m}+h(f_m(w_m), f(w_m))+ h(f(w_m),
f(x_0))\rightarrow 0\,,\qquad m\rightarrow\infty\,,
\end{equation}
оскільки $f_m$ збігається до $f$ рівномірно у $D,$ а $f$ неперервна
в точці $x_0$ за доведеним вище. З огляду на те, що $f_m(x_0)\in
\partial D$ при $m\in {\Bbb N}$ за умовою відкритості, дискретності
і замкненості $f_m$ (див. \cite[теорема~3.3]{Vu}), а також по
замкненості межі $\partial D,$ маємо: $f(x_0)\in \partial
D^{\,\prime}.$ Тоді по неперервності $f$ у $\overline{D}$ маємо:
$f(x_m)=z_m\rightarrow z_0\in D^{\,\prime}$ і $f(x_0)=z_0\in
\partial D^{\,\prime}.$ Це суперечить обранню $z_0\in E\subset
D^{\,\prime}.$ Остаточно, $f(E)$ є замкненим в $D^{\,\prime},$ бо в
випадку $x_0\in D$ це встановлено, а випадок $x_0\in\partial D$
неможливий.

\medskip
Покажемо, що $f(\partial D)\subset \partial D^{\,\prime}.$ Справді,
нехай $x_0\in
\partial D.$ Треба довести, що $f(x_0)\in
\partial D^{\,\prime}.$ Оскільки $f_m$ має неперервне продовження в
точку $x_0,$ знайдеться послідовність $w_m\in D,$ $w_m\rightarrow
x_0$ при $m\rightarrow\infty,$ така що $h(f_m(w_m),
f_m(x_0))<\frac{1}{m}.$ Повторюючи міркування у~(\ref{eq6}),
приходимо до висновку $f(x_0)\in \partial D^{\,\prime},$ що і треба
було встановити.

\medskip
Покажемо, що $\overline{f}(\overline{D})=\overline{D^{\,\prime}}.$
Доведення аналогічно до останньої частини доведення теореми~3.1
в~\cite{SSD}. Очевидно, що
$\overline{f}(\overline{D})\subset\overline{D^{\,\prime}}.$
Покажемо, що $\overline{D^{\,\prime}}\subset
\overline{f}(\overline{D}).$ Справді, нехай $y_0\in
\overline{D^{\,\prime}},$ тоді або $y_0\in D^{\,\prime},$ або
$y_0\in \partial{D^{\,\prime}}.$ а) Якщо $y_0\in D^{\,\prime},$ то
$y_0=f_j(x_j),$ $x_j\in D.$ Можна вважати, що $x_j\rightarrow x_0,$
$j\rightarrow\infty,$ $x_0\in \overline{D}.$ За неперервністю $f$ у
$\overline{D}$ маємо: $f(x_0)=\lim\limits_{j\rightarrow\infty}
f(x_j).$ Проте, за рівномірною збіжністю $f_j$ до $f$ маємо:
$f_j(x_j)-f(x_j)=y_0-f(x_j)\rightarrow 0,$ $j\rightarrow\infty.$
Тоді $y_0=f(x_0).$ Оскільки за доведеним $f(\partial D)\subset
\partial{D^{\,\prime}},$ маємо: $x_0\in D.$ Отже, $y_0=f(x_0),$ $x_0\in
D.$ б) Нарешті, нехай $y_0\in
\partial D^{\,\prime},$ тоді знайдеться послідовність $y_k\in
D^{\,\prime}$ така, що $y_k=f(x_k)\rightarrow y_0$ при
$k\rightarrow\infty$ і $x_k\in D.$ Через компактність простору
$\overline{{\Bbb R}^n}$ ми можемо вважати, що $x_k\rightarrow x_0,$
де $x_0\in\overline{D}.$ Тоді $\overline{f}(x_0)=y_0\in
\overline{f}(\partial D)\subset \overline{f}(\overline{D}).$ Отже,
$\overline{D^{\,\prime}}\subset \overline{f}(\overline{D}).$

\medskip
Залишилося показати нульвимірність $\overline{f}$ на $\partial
D^{\,\prime}.$ Припустимо супротивне, а саме, нехай існує точка
$y_0\in
\partial D^{\,\prime}$ така, що $f^{\,-1}(y_0)\supset K_0,$ де
$K_0\subset\overline{D}$ -- деякий невироджений континуум. Оскільки
$f$ нульвимірне в $D,$ можна вважати, що $K_0\subset\partial D.$

\medskip
Оскільки $\overline{D}$ -- компакт в $\overline{{\Bbb R}^n},$ і крім
того, відображення $\overline{f}$ є неперервним в $\overline{D},$
воно є рівномірно неперервним в $\overline{D}.$ В цьому випадку, для
кожного $j\in {\Bbb N}$ знайдеться $\delta_j<1/j$ таке, що
\begin{equation}\label{eq3K}
h(\overline{f}(x),\overline{f}(x_0))=h(\overline{f}(x),y_0)<1/j
\quad \forall\,\, x,x_0\in \overline{D},\quad h(x, x_0)<\delta_j\,,
\quad \delta_j<1/j\,.
\end{equation}
Позначимо $B_h(x_0, r)=\{x\in \overline{{\Bbb R}^n}: h(x, x_0)<r\}.$
Тоді для заданого $j\in {\Bbb N},$ покладемо
$$B_j:=\bigcup\limits_{x_0\in K_0}B_h(x_0, \delta_j)\,,\quad j\in {\Bbb N}\,.$$
Оскільки множина $B_j$ є околом континуума~$K_0,$
за~\cite[лема~2.2]{HK} існує окіл $U_j$ множини $K_0$ такий, що
$U_j\subset B_j$ і множина $U_j\cap D$ є зв'язною. Без обмеження
загальності можна вважати, що $U_j$ -- відкрита множина. Тоді
множина $U_j\cap D$ також лінійно зв'язна
(див.~\cite[пропозиція~13.1]{MRSY}). Оскільки $K_0$ -- компакт, то
знайдуться $z_0, w_0\in K_0$ такі, що $h(K_0)=h(z_0, w_0).$ Звідси
випливає, що знайдуться елементи $z_j\in U_j\cap D$ і $w_j\in
U_j\cap D$ такі, що $z_j\rightarrow z_0$ і $w_j\rightarrow w_0$ при
$j\rightarrow\infty.$ Ми можемо припускати, що
\begin{equation}\label{eq2B}
h(z_j, w_j)>h(K_0)/2\quad \forall\,\, j\in {\Bbb N}\,.
\end{equation}
Оскільки множина $U_j\cap D$ лінійно зв'язна, точки $z_j$ і $w_j$
деякою кривою $\gamma_j\in U_j\cap D.$ Покладемо $C_j:=|\gamma_j|.$

\medskip
Зауважимо, що $h(f(C_j))\rightarrow 0$ при $j\rightarrow\infty.$
Дійсно, оскільки $f(C_j)$ -- континуум в $\overline{{\Bbb R}^n},$
існують такі точки $y_j, y^{\,\prime}_j\in f(C_j)$ такі, що
$h(f(C_j))=h(y_j, y^{\,\prime}_j).$ Тоді існують $x_j,
x^{\,\prime}_j\in C_j$ такі, що $y_j=f(x_j)$ і
$y^{\,\prime}_j=f(x^{\,\prime}_j).$ Тоді точки $x_j$ і
$x^{\,\prime}_j$ належать до $U_j\subset B_j.$ Це означає, що
існують $x^j_1$ і $x^j_2\in K_0$ такі, що $x_j\in B(x^j_1,
\delta_j)$ і $x^{\,\prime}_j\in B(x^j_2, \delta_j).$ В такому
випадку, за співвідношенням~(\ref{eq3K}) і з огляду на нерівність
трикутника ми отримаємо, що
$$h(f(C_j))=h(y_j, y^{\,\prime}_j)=h(f(x_j), f(x^{\,\prime}_j))\leqslant$$
\begin{equation}\label{eq3L}
\leqslant h(f(x_j), f(x^j_1))+h(f(x^j_1), f(x^j_2))+h(f(x^j_2),
f(x^{\,\prime}_j))<2/j\rightarrow 0\quad {\text при}\quad
j\rightarrow\infty\,.\end{equation}
Зі співвідношень~(\ref{eq2B}) і~(\ref{eq3L}) випливає, що континууми
$C_j,$ $j=1,2,\ldots ,$ задовольняють умови леми~\ref{lem1}. За цією
лемою $h(f(C_j), y_0)\geqslant \delta_1>0$ для всіх $j\in {\Bbb N}.$
З іншого боку, оскільки за доведеним вище $x_j\in B(x^j_1,
\delta_j),$ то з огляду на співвідношення~(\ref{eq3K}) ми отримаємо,
що $h(f(x_j),y_0)<1/j,$ $j=1,2,\ldots .$ Отримана суперечність
вказує на невірність припущення щодо відсутності нульвимірності
відображення $\overline{f}$ у $\overline{D}.$ Лема доведена.~$\Box$
\end{proof}

\medskip
Нарешті, сформулюємо і доведемо ключове твердження про дискретність
відображення (див. \cite[теорема~4.7]{Vu}).

\medskip
\begin{lemma}\label{lem3}
{\sl\,Нехай за умов леми~\ref{lem2} область $D^{\,\prime}$ є
локально зв'язною на своїй межі, а відображення $f$ зберігає
орієнтацію.  Тоді відображення $f$ є відкритим, дискретним і
замкненим у $D,$ крім того, $f$ має неперервне продовження
$\overline{f}:\overline{D}\rightarrow \overline{D^{\,\prime}},$
$\overline{f}(\overline{D})=\overline{D^{\,\prime}}$ таке, що $N(f,
D)=N(f, \overline{D})<\infty.$ Зокрема, $\overline{f}$ є дискретним
у $\overline{D}.$ }
\end{lemma}

\medskip
\begin{proof}
Доведення майже повністю повторює схему, викладену
в~\cite[теорема~4.7]{Vu}, див. також теорему~6.1 в~\cite{Sev$_2$}.
Передусім, можливість продовження $f$ до неперервного відображення
$\overline{f}:\overline{D}\rightarrow \overline{D^{\,\prime}},$
нульвимірність і замкненість $f$ і рівність
$\overline{f}(\overline{D})=\overline{D^{\,\prime}}$ випливають з
леми~\ref{lem2}. Оскільки $f$ зберігає орієнтацію і нульвимірне,
воно також відкрите і дискретне за теоремою Тітуса-Янга (див.
\cite[стор.~333]{TY}). Зауважимо також, що $N(f, D)<\infty$ з огляду
на~\cite[теорема~2.8]{MS}. Залишилося довести рівність $N(f, D)=N(f,
\overline{D}).$ Далі будемо міркувати методом від супротивного:
припустимо, що вказана рівність не виконується. Тоді знайдуться
точки $y_0\in
\partial D^{\,\prime}$ і $x_1,x_2,\ldots, x_k, x_{k+1}\in \partial
D$ такі, що $f(x_i)=y_0,$ $i=1,2,\ldots, k+1$ і $k:=N(f, D).$ Будемо
вважати, що $y_0\ne\infty.$ Оскільки область $D^{\,\prime}$ є
локально зв'язною в усіх точках своєї межі за припущенням, то для
кожного $p\in {\Bbb N}$ знайдеться окіл
$\widetilde{U^{\,\prime}_p}\subset B(y_0, 1/p)$ такий, що множина
$\widetilde{U^{\,\prime}_p}\cap D^{\,\prime}=U^{\,\prime}_p$ є
зв'язною.

Доведемо, що для кожного $i=1,2,\ldots, k+1$ знайдеться компонента
$V_p^i$ множини $f^{\,-1}(U^{\,\prime}_p)$ така, що
$x_i\in\overline{V_p^i}.$ Зафіксуємо $i=1,2,\ldots, k+1.$ За
неперервністю $\overline{f}$ у $\overline{D}$ знайдеться
$r_i=r_i(x_i)>0$ таке, що $f(B(x_i, r_i)\cap D)\subset
U^{\,\prime}_p.$ Оскільки область, яка є слабко плоскою, є локально
зв'язною в кожній точці своїй межі (див.~\cite[лема~3.15]{MRSY}), то
знайдеться окіл $W_i\subset B(x_i, r_i)$ точки $x_i$ такий, що
$W_i\cap D$ зв'язно. Тоді $W_i\cap D$ належить одній і тільки одній
компоненті зв'язності $V^p_i$ множини $f^{\,-1}(U^{\,\prime}_p),$
причому $x_i\in\overline{W_i\cap D}\subset \overline{V_p^i},$ що і
треба було довести.

Далі покажемо, що множини $\overline{V_p^i}$ є непересічними при
всіх $i=1,2,\ldots, k+1$ і достатньо великих $p\in {\Bbb N}.$ В свою
чергу, доведемо для цього, що $h(\overline{V_p^i})\rightarrow 0$ при
$p\rightarrow\infty$ для кожного фіксованого $i=1,2,\ldots, k+1.$
Проведемо доведення від супротивного. Тоді існує таке $1\leqslant
i_0\leqslant k+1,$ число $r_0>0,$
$r_0<\frac{1}{2}\min\limits_{1\leqslant i, j\leqslant k+1, i\ne
j}h(x_i, x_j)$ і зростаюча послідовність номерів $p_m,$
$m=1,2,\ldots,$ такі, що $S_h(x_{i_0}, r_0)\cap
\overline{V_{p_m}^{i_0}}\ne\varnothing,$ де $S_h(x_0, r)=\{x\in
\overline{{\Bbb R}^n}: h(x, x_0)=r\}.$ В такому випадку, існують
$a_m, b_m\in V_{p_m}^{i_0},$ такі що $a_m\rightarrow x_{i_0}$ при
$m\rightarrow\infty$ і $h(a_m, b_m)\geqslant r_0/2.$ З'єднаємо точки
$a_m$ і $b_m$ кривою $C_m,$ яка цілком належить $V_{p_m}^{i_0}.$
Тоді $h(|C_m|)\geqslant r_0/2$ при $m=1,2,\ldots .$ З іншого боку,
оскільки $|C_m|\subset f(V_{p_m}^{i_0})\subset B(y_0, 1/p_m),$ то
одночасно $h(f(|C_m|))\rightarrow 0$ при $m\rightarrow\infty$ і
$h(|C_m|, y_0)\rightarrow 0$ при $m\rightarrow\infty,$ що суперечить
лемі~\ref{lem1}. Отримана суперечність вказує на невірність
зробленого вище припущення.

За~\cite[лема~3.6]{Vu} $f$ є відображенням кожного
$\overline{V_p^i},$ $i=1,2,\ldots, k, k+1,$ на все $U^{\,\prime}_p,$
тобто, $N(f, D)\geqslant k+1,$ що суперечить означення числа $k.$
Отримана суперечність спростовує припущення, що $N(f,
\overline{D})>N(f, D).$ Лема доведена.~$\Box$
\end{proof}

\medskip
Доведення теорем~\ref{th1} і~\ref{th3} випливає з лем~\ref{lem2}
і~\ref{lem3} за рахунок спеціально підібраної функції $\psi,$ див.
\cite[лема~1.3]{Sev$_3$}.~$\Box$

\medskip
{\bf 3. Випадок простих кінців.}  Результати, доведені в попередніх
секціях, залишаються справедливими для відображених областей, які не
є локально зв'язними на своїй межі, але за певної додаткової
геометрії на ці області. Для того, щоб сформулювати ці результати,
наведемо деякі означення.

\medskip
Нехай $\omega$ --  відкрита множина в ${\Bbb R}^k$,
$k=1,\ldots,n-1$. Неперервне відображення
$\sigma\colon\omega\rightarrow{\Bbb R}^n$ називається {\it
$k$-вимірною поверхнею} в ${\Bbb R}^n$. {\it Поверхнею} будемо
називати довільну $(n-1)$-вимірну поверхню $\sigma$ в ${\Bbb R}^n.$
Поверхня $\sigma$ називається {\it жордановою поверхнею}, якщо
$\sigma(x)\ne\sigma(y)$ при $x\ne y$. Далі ми іноді будемо
використовувати $\sigma$ для позначення всього образу
$\sigma(\omega)\subset {\Bbb R}^n$ при відображенні $\sigma$,
$\overline{\sigma}$ замість $\overline{\sigma(\omega)}$ в ${\Bbb
R}^n$ і $\partial\sigma$ замість
$\overline{\sigma(\omega)}\setminus\sigma(\omega)$. Жорданова
поверхня $\sigma\colon\omega\rightarrow D$ в області $D$ називається
{\it розрізом} області $D$, якщо $\sigma$ розділяє $D$, тобто
$D\setminus \sigma$ має більше однієї компоненти,
$\partial\sigma\cap D=\varnothing$ і $\partial\sigma\cap\partial
D\ne\varnothing$.

Послідовність $\sigma_1,\sigma_2,\ldots,\sigma_m,\ldots$ розрізів
області $D$ називається {\it ланцюгом}, якщо:

(i) множина $\sigma_{m+1}$ міститься в точності в одній компоненті
$d_m$ множини $D\setminus \sigma_m$, при цьому, $\sigma_{m-1}\subset
D\setminus (\sigma_m\cup d_m)$; (ii)
$\bigcap\limits_{m=1}^{\infty}\,d_m=\varnothing$. З означення
ланцюгу розрізів випливає, що $d_1\supset d_2\supset
d_3\supset\ldots\supset d_{m-1}\supset d_m\supset
d_{m+1}\supset\ldots\,.$
Два ланцюги розрізів $\{\sigma_m\}$ і $\{\sigma_k^{\,\prime}\}$
називаються {\it еквівалентними}, якщо для кожного $m=1,2,\ldots$
область $d_m$ містить всі області $d_k^{\,\prime}$ за виключенням
скінченної кількості, і для кожного $k=1,2,\ldots$ область
$d_k^{\,\prime}$ також містить всі області $d_m$ за виключенням
скінченної кількості.

{\it Кінець} області $D$ --- це клас еквівалентних ланцюгів розрізів
області $D$. Нехай $K$ --- кінець області $D$ в ${\Bbb R}^n$, тоді
множина $I(K)=\bigcap\limits_{m=1}\limits^{\infty}\overline{d_m}$
називається {\it тілом кінця} $K$. Скрізь далі, як зазвичай,
$\Gamma(E, F, D)$ позначає сім'ю всіх таких кривих $\gamma\colon[a,
b]\rightarrow D$, що $\gamma(a)\in E$ і $\gamma(b)\in F,$ крім того,
$M(\Gamma)$ позначає модуль сім'ї кривих $\Gamma$ в ${\Bbb R}^n,$ а
запис $\rho\in {\rm adm}\,\Gamma$ означає, що функція $\rho$
борелева, невід'ємна і має довжину, не меншу одиниці в метриці
$\rho$ (див.~\cite{Na}, \cite{Va}). Слідуючи~\cite{Na}, будемо
говорити, що кінець $K$ є {\it простим кінцем}, якщо $K$ містить
ланцюг розрізів $\{\sigma_m\}$, такий, що $M(\Gamma(\sigma_m,
\sigma_{m+1}, D))<\infty$ при всіх $m\in {\Bbb N}$ і $
\lim\limits_{m\rightarrow\infty}M(\Gamma(C, \sigma_m, D))=0 $
для деякого континууму $C$ в $D.$ 
%
Далі використовуються наступні позначення: множина простих кінців,
що відповідають області $D,$ позначається символом $E_D,$ а
поповнення області $D$ її простими кінцями позначається
$\overline{D}_P.$ Будемо говорити, що межа області $D$ в ${\Bbb
R}^n$ є {\it локально квазіконформною}, якщо кожна точка
$x_0\in\partial D$ має окіл $U$ в ${\Bbb R}^n$, який може бути
відображений квазіконформним відображенням $\varphi$ на одиничну
кулю ${\Bbb B}^n:=B(0, 1)\subset{\Bbb R}^n$ так, що
$\varphi(\partial D\cap U)$ є перетином ${\Bbb B}^n$ з координатною
гіперплощиною. Для множин $E\subset {\Bbb R}^n$ і $A, B\subset {\Bbb
R}^n$ покладемо
$$d(E):=\sup\limits_{x, y\in E}|x-y|\,,\quad d(A, B):=
\inf\limits_{x\in A, y\in B}|x-y|\,.$$
Будемо називати ланцюг розрізів $\{\sigma_m\}$ {\it регулярним},
якщо $\overline{\sigma_m}\cap\overline{\sigma_{m+1}}=\varnothing$
при кожному $m\in {\Bbb N}$ і, крім того, $d(\sigma_{m})\rightarrow
0$ при $m\rightarrow\infty.$ Якщо кінець $K$ містить принаймні один
регулярний ланцюг, то $K$ будемо називати {\it регулярним}.
Говоримо, що обмежена область $D$ в ${\Bbb R}^n$ {\it регулярна},
якщо $D$ може бути квазіконформно відображена на область з локально
квазіконформною межею, замикання якої є компактом в ${\Bbb R}^n,$
крім того, кожен простий кінець $P\subset E_D$ є регулярним.
Зауважимо, що у просторі ${\Bbb R}^n$ кожний простий кінець
регулярної області містить ланцюг розрізів з властивістю
$d(\sigma_{m})\rightarrow 0$ при $m\rightarrow\infty,$ і навпаки,
якщо у кінця є вказана властивість, то він -- простий
(див.~\cite[лема~3.1]{IS}; див. також~\cite[теорема~5.1]{Na}). Крім
того, замикання $\overline{D}_P$ регулярної області $D$ є {\it
метризовним}, при цьому, якщо $g:D_0\rightarrow D$ -- квазіконформне
відображення області $D_0$ з локально квазіконформною межею на
область $D,$ то для $x, y\in \overline{D}_P$ покладемо:
\begin{equation}\label{eq1A}
\rho(x, y):=|g^{\,-1}(x)-g^{\,-1}(y)|\,,
\end{equation}
де для $x\in E_D$ елемент $g^{\,-1}(x)$ розуміється як деяка (єдина)
точка межі $D_0,$ коректно визначена з огляду
на~\cite[теорема~2.1]{IS}; див. також~\cite[теорема~4.1]{Na}.
Зокрема, будемо говорити, що послідовність $x_m\in D,$
$m=1,2,\ldots,$ {\it збігається} до простого кінця $P\in E_D$ при
$m\rightarrow\infty,$ якщо для будь-якого натурального $k\in {\Bbb
N}$ всі елементи послідовності $x_m,$ крім скінченної кількості,
належать області $d_k$ (де $d_k,$ $k=1,2,\ldots$ -- послідовність
вкладених областей з означення простого кінця $P$).

\medskip
Справедлива наступна

\medskip
\begin{lemma}\label{lem4}
{\sl\, Нехай $D, D^{\,\prime}\subset {\Bbb R}^n,$ $n\geqslant 2,$ --
області в ${\Bbb R}^n,$ причому $D$ має слабко плоску межу, а
$D^{\,\prime}$ є регулярною. Припустимо, $f_m,$ $m=1,2,\ldots ,$ --
послідовність відкритих дискретних відображень області $D$ на
$D^{\,\prime},$ які збігаються рівномірно у $D$ до деякого (не
тотожно сталого) відображення $f:D\rightarrow
\overline{D^{\,\prime}}_P$ по метриці $\rho$ у
$\overline{D^{\,\prime}}_P,$ тобто,
$$\sup\limits_{x\in D}\rho(f_j(x), f(x))\rightarrow 0\,,\qquad j\rightarrow
\infty\,,$$
і задовольняють співвідношення~(\ref{eq2*A}) в кожній точці $y_0\in
\overline{D^{\,\prime}}.$ Припустимо, що для кожного $y_0\in
\overline{D^{\,\prime}}$ знайдеться
$\varepsilon_0=\varepsilon_0(y_0)>0$ і вимірна за Лебегом функція
$\psi:(0, \varepsilon_0)\rightarrow [0,\infty]$ такі, що виконані
умови~(\ref{eq7***})--(\ref{eq3.7.2}), де $A(y_0, \varepsilon,
\varepsilon_0)$ визначено в (\ref{eq1**}). Нехай $C_j,$
$j=1,2,\ldots ,$ -- довільна послідовність континуумів така, що
$h(C_j)\geqslant \delta>0$ для деякого $\delta>0$ і всіх $j\in {\Bbb
N}$ і, крім того, $\rho(f(C_j))\rightarrow 0$ при
$j\rightarrow\infty.$ Тоді для кожного $P_0\in
\overline{D^{\,\prime}}_P$ знайдеться $\delta_1>0$ таке, що
\begin{equation}\label{eq7}
\rho(f(C_j), P_0)\geqslant \delta_1>0
\end{equation}
для всіх $m\in {\Bbb N}$ і $j\in {\Bbb N}.$ Тут, як звично,
$$\rho(A)=\sup\limits_{x,
y\in A}\rho(x, y)\,,\qquad \rho(A, B)=\inf\limits_{x\in A, y\in
B}\rho(x, y)\,.$$ }
\end{lemma}
\begin{proof}
Надалі вважаємо $P_0\ne\infty.$ Припустимо супротивне, а саме, нехай
для деякої зростаючої послідовності номерів $j_k,$ $k=1,2,\ldots ,$
виконана умова: $\rho(f(C_{j_k}), P_0)\rightarrow 0$ при
$k\rightarrow\infty.$ Оскільки за припущенням $D^{\,\prime}$
регулярна, існує відображення $g$ деякої області $D_0$ з локально
квазіконформною межею на $D^{\,\prime}.$

\medskip
Доведемо, що зі збіжності послідовності $f_j$ до відображення $f$
рівномірно по метриці $\rho$ при $j\rightarrow\infty$ випливає
локально рівномірна збіжність цієї ж послідовності у $D$ відносно
звичайної евклідової метрики. Припустимо протилежне. Тоді знайдеться
компакт $C$ у $D,$ на якому збіжність $f_j$ до $f$ не рівномірна.
Звідси випливає наявність числа $\varepsilon_0>0,$ точок $x_k\in C,$
$k=1,2,\ldots,$ зростаючої послідовності номерів $j_k,$
$k=1,2,\ldots,$ і відображень $f_{j_k},$ $k=1,2,\ldots,$ таких, що
\begin{equation}\label{eq1B}
|f_{j_k}(x_k)-f(x_k)|\geqslant\varepsilon_0\,.
\end{equation}
Оскільки $f_j$ збігається до $f$ рівномірно в метриці $\rho,$ це
означає, що
\begin{equation}\label{eq1C}
g^{\,-1}(f_{j_k}(x_k))-g^{\,-1}(f(x_k))\rightarrow 0\,, \qquad
k\rightarrow\infty\,.
\end{equation}
Оскільки $C$ -- компакт, то $g^{\,-1}(f(C))$ також компакт у
$D^{\,\prime}$ як неперервний образ компакту. Тоді існує збіжна
підпослідовність послідовності $g^{\,-1}(f(x_k)),$ $k=1,2,\ldots .$
Без обмеження загальності можна вважати, що сама ця послідовність
збігається до деякої точки $\omega_0\in D^{\,\prime}$ при
$k\rightarrow\infty.$ Тоді зі співвідношення~(\ref{eq1C}) випливає,
що
\begin{equation}\label{eq1D}
g^{\,-1}(f_{j_k}(x_k))-\omega_0\rightarrow 0\,, \qquad
k\rightarrow\infty\,.
\end{equation}
Оскільки $g$ -- гомеомоморфізм у $D^{\,\prime},$ з~(\ref{eq1D})
випливає, що $f_{j_k}(x_k)\rightarrow g(\omega_0).$ З іншого боку,
співвідношення $g^{\,-1}(f(x_k))\rightarrow \omega_0$ при
$k\rightarrow \infty$ тягне за собою, що $f(x_k)\rightarrow
g(\omega_0).$ За нерівністю трикутника
$$|f_{j_k}(x_k)-f(x_k)|\leqslant |f_{j_k}(x_k)-g(\omega_0)|
+|g(\omega_0)-f(x_k)|\rightarrow 0$$
при $k\rightarrow\infty.$ Останнє співвідношення
суперечить~(\ref{eq1B}), що і доводить локально рівномірну збіжність
сім'ї $f_j,$ $j=1,2,\ldots ,$ у $D.$

\medskip
Тоді з огляду на~\cite[лема~1]{SevSkv$_3$} $f$ або нульвимірне, або
стале відображення. Оскільки за умовою $f$ не є сталим, $f$ --
нульвимірне. Тоді існують точка $z_0\in D,$ число $r_0>0$ и окіл
елемента $P_0$ такі, що
\begin{equation}\label{eq3I}
\overline{B(z_0, r_0)}\subset D\,,\qquad f(\overline{B(z_0,
r_0)})\cap \overline{U}=\varnothing\,.
\end{equation}
Покладемо $F:=\overline{B(z_0, r_0)}.$ З огляду на~(\ref{eq3I})
покажемо, що існує окіл $V$ елементу $P_0$ такий, що
\begin{equation}\label{eq5}
f_j(F)\cap V=\varnothing\,,\qquad j=1,2,\ldots\,.
\end{equation}
Припустимо протилежне. Тоді існують послідовності $j_m\in{\Bbb N}$ і
$x_m\in F,$ $m=1,2,\ldots, $ такі, що $f_{j_m}(x_m)\rightarrow P_0$
при $m\rightarrow\infty.$ З огляду на те, що $f_{j_m}(F)$ є
компактом у $D^{\,\prime}$ при кожному фіксованому $m,$
послідовність індексів $j_m,$ $m=1,2,\ldots,$ можна вважати
зростаючою. Тоді за нерівністю трикутника
\begin{equation}\label{eq6A}
\rho(f(x_m), P_0)\leqslant \rho(f(x_m), f_{j_m}(x_m))+
\rho(f_{j_m}(x_m), P_0)\rightarrow 0
\end{equation}
при $m\rightarrow\infty,$ оскільки $\rho(f(x_m), f_{j_m}(x_m))$ при
$m\rightarrow\infty$ по рівномірній збіжності $f_m$ до $f$ у $D$ по
метриці $\rho.$ Співвідношення~(\ref{eq6A}) суперечить~(\ref{eq3I}),
що доводить~(\ref{eq5}).

\medskip
Позначимо $\Gamma_k:=\Gamma(F, C_{j_k}, D).$ Зауважимо, що область
зі слабко плоскою межею є рівномірною (див.
\cite[наслідок~4.3]{Sev$_2$}), тобто, якщо континууми $C_{j_k}$
задовольняють умову $h(C_{j_k})\geqslant \delta>0$ і
$h(F)>\delta_*>0,$ то знайдеться $\delta_2>0$ таке, що
\begin{equation}\label{eq3P}
M(\Gamma_k)\geqslant \delta_2>0
\end{equation}
для всіх $k\in {\Bbb N}.$ Нехай $d_l,$ $l=1,2,\ldots ,$ --
послідовність областей, яка відповідає простому кінцю $P_0,$ і нехай
$\sigma_l$ -- розріз, що відповідає $d_l.$ Можна вважати, що усі
$\sigma_l,$ $l=1,2,\ldots, $ лежать на сферах $S(y_0, r_l )$ з
центром у деякій точці $y_0\in \partial D^{\,\prime},$ де
$r_l\rightarrow 0$ при $l\rightarrow\infty$
(див.~\cite[лема~3.1]{IS}, див. також~\cite[лема~1]{KR}).

\medskip Доведемо, що для кожного $l\in {\Bbb
N}$ існує номер $k=k_l$ такий, що
\begin{equation}\label{eq3Q}
f_{j_k}(C_{j_k})\subset d_l\,,\qquad k\geqslant k_l\,.
\end{equation}
Припустимо протилежне. Тоді існує $l_0\in {\Bbb N}$ таке, що
\begin{equation}\label{eq3R}
f_{j_{k_m}}(C_{j_{k_m}})\cap ({\Bbb R}^n\setminus
d_{l_0})\ne\varnothing
\end{equation}
для деякої зростаючої послідовності номерів $k_m,$ $m=1,2,\ldots .$
В такому випадку, знайдеться послідовність $z_{m}\in
f_{j_{k_m}}(C_{j_{k_m}})\cap ({\Bbb R}^n\setminus d_{l_0}),$ $m\in
{\Bbb N}.$ Нехай $z_m=f_{j_{k_m}}(x_m),$ $x_m\in C_{j_{k_m}}.$

Оскільки за припущенням $\rho(f(C_{j_k}), P_0)\rightarrow 0$ для
деякої зростаючої послідовності номерів $j_k,$ $k=1,2,\ldots ,$ то
зокрема
\begin{equation}\label{eq3S}
\rho(f(C_{j_{k_m}}), P_0)\rightarrow 0\qquad {\text при}\qquad
m\rightarrow\infty\,.
\end{equation}
Оскільки $\rho(f(C_{j_{k_m}}), P_0)=\inf\limits_{y\in
f(C_{j_{k_m}})}\rho(y, P_0)$ і $f(C_{j_{k_m}})$ є компактом як
неперервний образ компакту $C_{j_{k_m}}$ при відображенні $f,$
звідси випливає, що $\rho(f(C_{j_{k_m}}), P_0)=\rho(y_m, P_0),$ де
$y_m\in f(C_{j_{k_m}}).$ Тоді зі співвідношення~(\ref{eq3S}) ми
отримаємо, що $y_m\rightarrow P_0$ при $m\rightarrow\infty.$ Нехай
$y_m=f(\xi_m),$ де $\xi_m\in C_{j_{k_m}}.$ Отже, за нерівністю
трикутника
$$\rho(z_m, P_0)=\rho(f_{j_{k_m}}(x_m), P_0)\leqslant$$
\begin{equation}\label{eq5E}
\leqslant \rho(f_{j_{k_m}}(x_m), f(x_m))+\rho(f(x_m),
f(\xi_m))+\rho(f(\xi_m), P_0)\rightarrow 0\,, \quad
m\rightarrow\infty\,,
\end{equation}
оскільки $\rho(f_{j_{k_m}}(x_m), f(x_m))\rightarrow 0$ при
$m\rightarrow\infty$ по рівномірній збіжності $f_m$ у $D$ по метриці
$\rho;$ крім того, $\rho(f(x_m), f(\xi_m))\leqslant
\rho(f(C_{j_{k_m}}))\rightarrow 0$ при $m\rightarrow\infty$ за умови
леми; нарешті, $\rho(f(\xi_m), P_0)=\rho(y_m, P_0)\rightarrow 0$ при
$m\rightarrow\infty$ за доведеним вище. Отже, з огляду
на~(\ref{eq5E}) $z_m\rightarrow P_0$ при $m\rightarrow\infty,$ але
це суперечить умові $z_{m}\in f_{j_{k_m}}(C_{j_{k_m}})\cap ({\Bbb
R}^n\setminus d_{l_0}).$ Ця умова, в свою чергу,
спростовує~(\ref{eq3R}) і доводить~(\ref{eq3Q}).

\medskip
З~(\ref{eq3Q}) випливає, що $f_{j_{k_l}}(C_{j_{k_l}})\subset d_l\,,$
причому можна вважати послідовність номерів $k_l,$ $l=1,2,\ldots ,$
зростаючою. Переходячи к перенумерації в разі потреби, не обмежуючи
загальності множна вважати, шо сама послідовність $k$ задовольняє ці
умови, тобто,
\begin{equation}\label{eq3T}
f_{j_{k}}(C_{j_{k}})\subset d_k\,, \qquad k=1,2,\ldots\,.
\end{equation}
З огляду на~(\ref{eq5}) можна вважати, що
\begin{equation}\label{eq3U}
f_{j_{k}}(F)\subset D^{\,\prime}\setminus d_1\,,\qquad
k=1,2,\ldots\,.
\end{equation}
В такому випадку, зауважимо, що
\begin{equation}\label{eq3GA}
f_{j_{k}}(\Gamma_{j_{k}})>\Gamma(S(y_0, r_k), S(y_0, r_1), A(y_0,
r_k, r_1))\,.
\end{equation}
Дійсно, нехай $\widetilde{\gamma}\in f_{j_{k}}(\Gamma_{j_{k}}).$
Тоді $\widetilde{\gamma}(t)=f_{j_{k}}(\gamma(t)),$ де $\gamma\in
\Gamma_{j_{k}},$ $\gamma:[0, 1]\rightarrow D,$ $\gamma(0)\in F,$
$\gamma(1)\in C_{j_{k}}.$ За співвідношенням~(\ref{eq3U}) маємо, що
$f_{j_{k}}(\gamma(0))\in f_{j_{k}}(F)\subset D^{\,\prime}\setminus
d_1,$ крім того, за співвідношенням~(\ref{eq3T}) ми отримаємо, що
$f_{j_{k}}(\gamma(1))\in f_{j_{k}}(C_{j_{k}})\subset d_k.$ Отже,
$|f(\gamma(t))|\cap d_1\ne\varnothing \ne |f(\gamma(t))|\cap
(D^{\,\prime}\setminus d_1).$ Тоді за \cite[теорема~1.I.5.46]{Ku} та
означенням розрізу області існує $0<t_1<1$ таке, що
$f_{j_{k}}(\gamma(t_1))\in
\partial d_1\cap D^{\,\prime}\subset S(y_0, r_1).$
Покладемо $\gamma_1:=\gamma|_{[t_1, 1]}.$ Можна вважати, що
$f_{j_{k}}(\gamma(t))\in d_1$ при всіх $t\geqslant t_1.$ Міркуючи
аналогічно, ми отримаємо точку $t_2\in (t_1, 1]$ таку, що
$f_{j_{k}}(\gamma(t_2))\in S(y_0, r_k).$ Покладемо
$\gamma_2:=\gamma|_{[t_1, t_2]}.$ Можна вважати, що
$f_{j_{k}}(\gamma(t))\in d_k$ при всіх $t\in [t_1, t_2].$ Тоді крива
$f_{j_{k}}(\gamma_2)$ є підкривою кривої
$f_{j_{k}}(\gamma)=\widetilde{\gamma},$ яка належить сім'ї
$\Gamma(S(y_0, r_k), S(y_0, r_1), A(y_0, r_k, r_1)).$
Співвідношення~(\ref{eq3GA}) встановлено.

\medskip
З~(\ref{eq3GA}) випливає, що
\begin{equation}\label{eq3HA}
\Gamma_{j_{k}}>\Gamma_{f_{j_{k}}}(S(y_0, r_k), S(y_0, r_1), A(y_0,
r_k, r_1))\,.
\end{equation}
Покладемо
$$\eta_{k}(t)=\left\{
\begin{array}{rr}
\psi(t)/I(r_k, r_1), & t\in (r_k, r_1)\,,\\
0,  &  t\not\in (r_k, r_1)\,,
\end{array}
\right. $$
де $I(r_k, r_1)=\int\limits_{r_k}^{r_1}\,\psi (t)\, dt.$ Зауважимо,
що
$\int\limits_{r_k}^{r_1}\eta_{l}(t)\,dt=1.$ Тоді за
співвідношеннями~(\ref{eq3.7.2}) і~(\ref{eq3HA})  з огляду на
означення відображення $f_{j_{k}}$ у~(\ref{eq2*A}) будемо мати:
$$M(\Gamma_{j_{k}})\leqslant M(\Gamma_{f_{j_{k}}}(S(y_0, r_k), S(y_0,
r_1), A(y_0,r_k, r_1)))\leqslant$$
\begin{equation}\label{eq3JA}
\leqslant \frac{1}{I^n(r_k, r_1)}\int\limits_{A(y_0, r_k, r_1)}
Q(y)\cdot\psi^{\,n}(|y-y_0|)\,dm(y)\rightarrow 0\quad
\text{при}\quad k\rightarrow\infty\,.
\end{equation}
Останнє співвідношення суперечить~(\ref{eq3A}). Отримана
суперечність доводить лему.~$\Box$
\end{proof}

\medskip
\begin{lemma}\label{lem5}
{\sl\, Нехай $D\subset {\Bbb R}^n,$ $n\geqslant 2,$ -- область, яка
має слабо плоску межу, а область $D^{\,\prime}\subset {\Bbb R}^n$ є
регулярною. Припустимо, $f_m,$ $m=1,2,\ldots ,$ -- послідовність
відкритих, дискретних і замкнених відображень області $D$ на
$D^{\,\prime},$ які збігаються рівномірно у $D$ до деякого
відображення $f:D\rightarrow {\Bbb R}^n$ і задовольняють
співвідношення~(\ref{eq2*A}) в кожній точці $y_0\in
\overline{D^{\,\prime}}.$ Нехай для будь-якого $y_0\in
\overline{D^{\,\prime}}$ знайдеться
$\varepsilon_0=\varepsilon_0(y_0)>0$ і вимірна за Лебегом функція
$\psi:(0, \varepsilon_0)\rightarrow [0,\infty]$ такі, що виконуються
співвідношення~(\ref{eq7***})--(\ref{eq3.7.2}), де $A(y_0,
\varepsilon, \varepsilon_0)$ визначено в (\ref{eq1**}). Тоді
відображення $f$ є нульвимірним і замкненим у $D,$ крім того, $f$
має неперервне продовження $\overline{f}:\overline{D}\rightarrow
\overline{D^{\,\prime}}_P,$
$\overline{f}(\overline{D})=\overline{D^{\,\prime}}_P$ і
$\overline{f}$ є нульвимірним відносно $\overline{D}_P.$}
\end{lemma}

\medskip
\begin{proof}
Можна вважати $y_0\ne \infty.$ Як було встановлено в ході доведення
леми~\ref{lem4}, послідовність $f_j$ збігається до $f$ локально
рівномірно в $D.$ Тоді з огляду на \cite{SevSkv$_3$} відображення
$f$ є нульвимірним у $D.$ За теоремою~1.1 і зауваженням~2.1
в~\cite{Sev$_4$}, кожне $f_m,$ $m=1,2,\ldots ,$ має неперервне
продовження $\overline{f_m}:\overline{D}\rightarrow
\overline{D^{\,\prime}}_P$ на $\overline{D}$ таке, що
$f_m(\overline{D})=\overline{D^{\,\prime}}_P.$ Покажемо, що існує
$\lim\limits_{x\rightarrow x_0}f(x)$ для кожного $x_0\in\partial D,$
де наявність границі слід розуміти в сенсі метрики $\rho$ у
$\overline{D^{\,\prime}}_P.$

Дійсно, за нерівністю трикутника
$$\rho(f(x), f(x_0))\leqslant \rho(f(x), f_m(x))+\rho(f_m(x), f_m(x_0))+\rho(f_m(x_0), f(x_0))\,.$$
У цій нерівності перший і третій доданки прямують до нуля з огляду
на рівномірну збіжність $f_m$ у $D.$ Отже, для будь-якого
$\varepsilon>0$ знайдеться номер $M=M(\varepsilon)$ такий, що при
$m\geqslant M(\varepsilon)$
\begin{equation}\label{eq5DA}
\rho(f(x), f(x_0))\leqslant \frac{2\varepsilon}{3}+\rho(f_m(x),
f_m(x_0))\,.
\end{equation}
При $m=M(\varepsilon)$ відображення $f_m$ неперервно продовжується в
точку $x_0,$ так що знайдеться $\delta=\delta(\varepsilon, x_0)>0$
таке, що $\rho(f_M(x), f_M(x_0))<\varepsilon/3$ при
$0<|x-x_0|<\delta.$ Тоді з огляду на~(\ref{eq5DA}) $\rho(f(x),
f(x_0))<\varepsilon$ при тих же $\delta.$

\medskip
Покажемо, що $f$ -- замкнене у $D.$ Нехай $E$ -- замкнене в $D.$
Треба показати, що $f(E)$ замкнене в $D^{\,\prime}.$ Нехай $z_m,$
$m=1,2,\ldots,$ -- послідовність у $f(E)$ така, що $z_m\rightarrow
z_0\in D^{\,\prime}$ при $m\rightarrow \infty.$ Покажемо, що $z_0\in
f(E).$ Оскільки $z_m,$ знайдуться $x_m\in E$ такі, що $f(x_m)=z_m,$
$m=1,2,\ldots .$ З огляду на компактність $\overline{{\Bbb R}^n}$
можна вважати, що послідовність $x_m$ збігається до деякого $x_0\in
\overline{D}$ при $m\rightarrow\infty.$ Якщо $x_0\in D,$ то $x_0\in
E,$ отже, $f(x_0)=z_0\in E$ з огляду на неперервність відображення
$f$ у $D.$ Якщо ж $x_0\in \partial D,$ то, оскільки $f_m$ має
неперервне продовження в точку $x_0,$ знайдеться послідовність
$w_m\in D,$ $w_m\rightarrow x_0$ при $m\rightarrow\infty,$ така що
$\rho(f_m(w_m), f_m(x_0))<\frac{1}{m}.$ Тоді за нерівністю
трикутника
$$\rho(f_m(x_0), f(x_0))\leqslant \rho(f_m(x_0), f_m(w_m))+\rho(f_m(w_m), f(w_m))+
\rho(f(w_m), f(x_0))\leqslant$$
\begin{equation}\label{eq6B}\leqslant\frac{1}{m}+\rho(f_m(w_m), f(w_m))+ \rho(f(w_m),
f(x_0))\rightarrow 0\,,\qquad m\rightarrow\infty\,,
\end{equation}
оскільки $f_m$ збігається до $f$ рівномірно у $D,$ а $f$ неперервна
в точці $x_0$ за доведеним вище. З огляду на те, що $f_m(x_0)\in
E_{D^{\,\prime}}$ при $m\in {\Bbb N}$ за умовою відкритості,
дискретності і замкненості $f_m$ (див. \cite[теорема~3.3]{Vu}), а
також по замкненості межі $\partial D,$ маємо: $f(x_0)\in
E_{D^{\,\prime}}.$ Тоді по неперервності $f$ у $\overline{D}$ маємо:
$f(x_m)=z_m\rightarrow z_0$ і $f(x_0)=z_0\in E_{D^{\,\prime}}.$ Це
суперечить обранню $z_0\in E\subset D^{\,\prime}.$ Остаточно, $f(E)$
є замкненим в $D^{\,\prime},$ бо в випадку $x_0\in D$ це
встановлено, а випадок $x_0\in\partial D$ неможливий.

\medskip
Покажемо, що $f(\partial D)\subset E_{D^{\,\prime}}.$ Нехай $x_0\in
\partial D.$ Треба довести, що $f(x_0)\in
E_{D^{\,\prime}}.$ Оскільки $f_m$ має неперервне продовження в точку
$x_0,$ знайдеться послідовність $w_m\in D,$ $w_m\rightarrow x_0$ при
$m\rightarrow\infty,$ така що $\rho(f_m(w_m),
f_m(x_0))<\frac{1}{m}.$ Повторюючи міркування у~(\ref{eq6B}),
приходимо до висновку $f(x_0)\in E_{D^{\,\prime}},$ що і треба було
встановити.

\medskip
Покажемо, що $\overline{f}(\overline{D})=\overline{D^{\,\prime}}_P.$
Очевидно, що
$\overline{f}(\overline{D})\subset\overline{D^{\,\prime}}_P.$
Покажемо, що $\overline{D^{\,\prime}}_P\subset
\overline{f}(\overline{D}).$ Справді, нехай $y_0\in
\overline{D^{\,\prime}}_P,$ тоді або $y_0\in D^{\,\prime},$ або
$y_0\in E_{D^{\,\prime}}.$ а) Якщо $y_0\in D^{\,\prime},$ то
$y_0=f_j(x_j),$ $x_j\in D.$ Можна вважати, що $x_j\rightarrow x_0,$
$j\rightarrow\infty,$ $x_0\in \overline{D}.$ За неперервністю $f$ у
$\overline{D}$ маємо: $f(x_0)=\lim\limits_{j\rightarrow\infty}
f(x_j).$ Проте, за рівномірною збіжністю $f_j$ до $f$ маємо:
$f_j(x_j)-f(x_j)=y_0-f(x_j)\rightarrow 0,$ $j\rightarrow\infty.$
Тоді $y_0=f(x_0).$ Оскільки за доведеним $f(\partial D)\subset
E_{D^{\,\prime}},$ маємо: $x_0\in D.$ Отже, $y_0=f(x_0),$ $x_0\in
D.$ б) Нарешті, нехай $y_0\in  E_{D^{\,\prime}},$ тоді знайдеться
послідовність $y_k\in D^{\,\prime}$ така, що $y_k=f(x_k)\rightarrow
y_0$ при $k\rightarrow\infty$ і $x_k\in D.$ Через компактність
простору $\overline{{\Bbb R}^n}$ ми можемо вважати, що
$x_k\rightarrow x_0,$ де $x_0\in\overline{D}.$ Тоді
$\overline{f}(x_0)=y_0\in \overline{f}(\partial D)\subset
\overline{f}(\overline{D}).$ Отже, $\overline{D^{\,\prime}}_P\subset
\overline{f}(\overline{D}).$

\medskip
Залишилося показати нульвимірність $\overline{f}$ на $\partial
D^{\,\prime}.$ Припустимо супротивне, а саме, нехай існує точка
$y_0\in
\partial D^{\,\prime}$ така, що $f^{\,-1}(y_0)\supset K_0,$ де
$K_0\subset\overline{D}$ -- деякий невироджений континуум. Тоді,
зокрема, $\overline{f}(K_0)=y_0.$ Оскільки за доведеним $f$ є
нульвимірним у $D,$ можна вважати, що $K_0\subset\partial D.$
Зауважимо, що $f\ne const$ з огляду на рівність
$\overline{f}(\overline{D})\subset\overline{D^{\,\prime}}_P.$

\medskip
Оскільки $\overline{D}$ -- компакт в $\overline{{\Bbb R}^n},$ і крім
того, відображення $\overline{f}$ є неперервним в $\overline{D},$
воно є рівномірно неперервним в $\overline{D}.$ В цьому випадку, для
кожного $j\in {\Bbb N}$ знайдеться $\delta_j<1/j$ таке, що
\begin{equation}\label{eq3KA}
\rho(\overline{f}(x),\overline{f}(x_0))=\rho(\overline{f}(x),y_0)<1/j
\quad \forall\,\, x,x_0\in \overline{D},\quad \rho(x,
x_0)<\delta_j\,, \quad \delta_j<1/j\,.
\end{equation}
Позначимо $B_h(x_0, r)=\{x\in \overline{{\Bbb R}^n}: h(x, x_0)<r\}.$
Тоді для заданого $j\in {\Bbb N},$ покладемо
$$B_j:=\bigcup\limits_{x_0\in K_0}B_h(x_0, \delta_j)\,,\quad j\in {\Bbb N}\,.$$
Оскільки множина $B_j$ є околом континуума~$K_0,$
за~\cite[лема~2.2]{HK} існує окіл $U_j$ множини $K_0$ такий, що
$U_j\subset B_j$ і множина $U_j\cap D$ є зв'язною. Без обмеження
загальності можна вважати, що $U_j$ -- відкрита множина. Тоді
множина $U_j\cap D$ також лінійно зв'язна
(see~\cite[Пропозиція~13.1]{MRSY}). Оскільки $K_0$ -- компакт, то
знайдуться $z_0, w_0\in K_0$ такі, що $h(K_0)=h(z_0, w_0).$ Звідси
випливає, що знайдуться елементи $z_j\in U_j\cap D$ і $w_j\in
U_j\cap D$ такі, що $z_j\rightarrow z_0$ і $w_j\rightarrow w_0$ при
$j\rightarrow\infty.$ Ми можемо припускати, що
\begin{equation}\label{eq2BA}
h(z_j, w_j)>h(K_0)/2\quad \forall\,\, j\in {\Bbb N}\,.
\end{equation}
Оскільки множина $U_j\cap D$ лінійно зв'язна, точки $z_j$ і $w_j$
деякою кривою $\gamma_j\in U_j\cap D.$ Покладемо $C_j:=|\gamma_j|.$

\medskip
Зауважимо, що $\rho(f(C_j))\rightarrow 0$ при $j\rightarrow\infty.$
Дійсно, оскільки $f(C_j)$ -- континуум в $\overline{{\Bbb R}^n},$
існують такі точки $y_j, y^{\,\prime}_j\in f(C_j)$ такі, що
$\rho(f(C_j))=\rho(y_j, y^{\,\prime}_j).$ Тоді існують $x_j,
x^{\,\prime}_j\in C_j$ такі, що $y_j=f(x_j)$ і
$y^{\,\prime}_j=f(x^{\,\prime}_j).$ Тоді точки $x_j$ і
$x^{\,\prime}_j$ належать до $U_j\subset B_j.$ Це означає, що
існують $x^j_1$ і $x^j_2\in K_0$ такі, що $x_j\in B(x^j_1,
\delta_j)$ і $x^{\,\prime}_j\in B(x^j_2, \delta_j).$ В такому
випадку, за співвідношенням~(\ref{eq3KA}) і з огляду на нерівність
трикутника ми отримаємо, що
$$\rho(f(C_j))=\rho(y_j, y^{\,\prime}_j)=\rho(f(x_j), f(x^{\,\prime}_j))\leqslant$$
\begin{equation}\label{eq3LA}
\leqslant \rho(f(x_j), f(x^j_1))+\rho(f(x^j_1),
f(x^j_2))+\rho(f(x^j_2), f(x^{\,\prime}_j))<2/j\rightarrow 0\quad
{\text при}\quad j\rightarrow\infty\,.\end{equation}
Зі співвідношень~(\ref{eq2BA}) і~(\ref{eq3LA}) випливає, що
континууми $C_j,$ $j=1,2,\ldots ,$ задовольняють умови
леми~\ref{lem1}. За цією лемою $\rho(f(C_j), y_0)\geqslant
\delta_1>0$ для всіх $j\in {\Bbb N}.$ З іншого боку, оскільки за
доведеним вище $x_j\in B(x^j_1, \delta_j),$ то з огляду на
співвідношення~(\ref{eq3KA}) ми отримаємо, що
$\rho(f(x_j),y_0)<1/j,$ $j=1,2,\ldots .$ Отримана суперечність
вказує на невірність припущення щодо відсутності нульвимірності
відображення $\overline{f}$ у $\overline{D}.$ Лема доведена.~$\Box$
\end{proof}

\medskip
\begin{lemma}\label{lem6}
{\sl\,Нехай за умов леми~\ref{lem5} відображення $f$ зберігає
орієнтацію. Тоді відображення $f$ є відкритим, дискретним і
замкненим у $D,$ крім того, $f$ має неперервне продовження
$\overline{f}:\overline{D}\rightarrow \overline{D^{\,\prime}}_P,$
$\overline{f}(\overline{D})=\overline{D^{\,\prime}}_P$ таке, що
$N(f, D)=N(f, \overline{D})<\infty.$ Зокрема, $\overline{f}$ є
дискретним у $\overline{D}.$ }
\end{lemma}

\medskip
\begin{proof}
Можливість неперервного продовження відображення $f$ до відображення
$\overline{f}:\overline{D}\rightarrow \overline{D^{\,\prime}}_P,$
його нульвимірність і замкненість, так само як і рівність
$\overline{f}(\overline{D})=\overline{D^{\,\prime}}_P$ є
твердженнями леми~\ref{lem5}. Оскільки $f$ зберігає орієнтацію і
нульвимірне, воно також відкрите і дискретне за теоремою Тітуса-Янга
(див. \cite[стор.~333]{TY}). Зауважимо, що $N(f, D) <\infty,$
див.~\cite[Теорема~2.8]{MS}. Доведемо що $N(f, D)=N(f,
\overline{D}).$ (Далі ми будемо міркувати, використовуючи схему
доведення теореми~4.7 у \cite{Vu}). Припустимо протилежне. Тоді
знайдуться точки $P_0\in E_{D^{\,\prime}}$ і $x_1,x_2,\ldots, x_k,
x_{k+1}\in
\partial D$ такі, що $f(x_i)=P_0,$ $i=1,2,\ldots, k+1$ і $k:=N(f,
D).$ Оскільки за припущенням $D^{\,\prime}$ є регулярною областю,
існує відображення $g$ області $D_0$ з локально квазіконформною
межею на $D^{\,\prime}.$ Розглянемо відображення $F:=f\circ
g^{\,-1}.$ Зауважимо, що за визначенням області з локально
квазіконформною межею $D_0$ є локально зв'язною на $\partial D_0.$ В
такому випадку, відображення $F$ має неперервне продовження
$\overline{F}:\overline{D}\rightarrow \overline{D_0}$ таке, що
$\overline{F}(\overline{D})=\overline{D_0}$ (це випливає з того, що
кожне з відображень $f$ і $g^{\,-1}$ має безперервне продовження на
$\overline{D}$ і $\overline{D^{\,\prime}}_P,$ відповідно). Покладемо
$y_0:=\overline{F}(P_0)\in D_0.$ Тоді для будь-якого $p\in {\Bbb N}$
існує окіл $\widetilde{U^{\,\prime}_p}\subset B(y_0, 1/p)$ точки
$y_0$ такий, що множина $\widetilde{U^{\,\prime}_p}\cap
D_0=U^{\,\prime}_p$ є зв'язною.

\medskip
Доведемо, що для будь-якого $i=1,2,\ldots, k+1$ існує компонента
$V_p^i$ множини $F^{\,-1}(U^{\,\prime}_p)$ така, що
$x_i\in\overline{V_p^i}.$ Зафіксуємо $i=1,2,\ldots, k+1.$ За
неперервністю $F$ в $\overline{D}$ існує $r_i=r_i(x_i)>0$ таке, що
$f(B(x_i, r_i)\cap D)\setminus U^{\,\prime}_p.$
За~\cite[лема~3.15]{MRSY} область зі слабко плоскою межею є локально
зв'язною на своїй межі. Отже, ми можемо знайти окіл $W_i\subset
B(x_i, r_i)$ точки $x_i$ такий, що $W_i\cap D$ є зв'язним. Тоді
$W_i\cap D$ належить одній і тільки одній компоненті $V^p_i$ множини
$F^{\,-1}(U^{\,\prime}_p),$ тоді як $x_i\in\overline{W_i\cap
D}\subset \overline{V_p^i},$ що і слід було встановити. Далі ми
покажемо, що множини $\overline{V_p^i}$ є непересічними для
будь-якого $i=1,2,\ldots, k+1$ і достатньо великому $p\in {\Bbb N}.$
У свою чергу, покажемо для цього, що $h(\overline{V_p^i})\rightarrow
0$ при $p\rightarrow\infty$ для кожного фіксованого $i=1,2,\ldots,
k+1.$ Припустимо протилежне. Тоді знайдеться $1\leqslant
i_0\leqslant k+1,$ і число $r_0>0,$
$r_0<\frac{1}{2}\min\limits_{1\leqslant i, j\leqslant k+1, i\ne
j}h(x_i, x_j)$ та зростаюча послідовність чисел $p_m,$
$m=1,2,\ldots,$ такі, що $S_h(x_{i_0}, r_0)\cap
\overline{V_{p_m}^{i_0}}\ne\varnothing,$ де $S_h(x_0, r)=\{x\in
\overline{{\Bbb R}^n}: h(x, x_0)=r\},$ і $h$ позначає хордальну
метрику в $\overline{{\Bbb R}^n}.$ У цьому випадку, знайдуться $a_m,
b_m\in V_{p_m}^{i_0}$ такі, що $a_m\rightarrow x_{i_0}$ при
$m\rightarrow\infty$ і $h(a_m, b_m)\geqslant r_0/2.$ З'єднаємо до
точки $a_m$ і $b_m$ кривою $C_m,$, яка повністю належить
$V_{p_m}^{i_0}.$ Тоді $h(|C_m|)\geqslant r_0/2$ при $m=1,2,\ldots .$
З іншого боку, оскільки $|C_m|\subset f(V_{p_m}^{i_0})\subset B(y_0,
1/p_m)$ і одночасно $h(F(|C_m|))\rightarrow 0$ при
$m\rightarrow\infty,$ $h(F(|C_m|), y_0)\rightarrow 0$ при
$m\rightarrow\infty.$ За визначенням метрики $\rho$ в~(\ref{eq5}) і
відображення $g,$ ми отримуємо що $\rho(f(|C_m|))\rightarrow 0$ при
$m\rightarrow\infty$ і $\rho(f(|C_m|), y_0)\rightarrow 0$ при
$m\rightarrow\infty,$ що суперечить лемі~\ref{lem4}. Отримане
протиріччя вказує на те, що $h(\overline{V_p^i})\rightarrow 0$ при
$p\rightarrow\infty$ для кожного фіксованого $i=1,2,\ldots, k+1.$ З
огляду на~\cite[лема~3.6]{Vu} $F$ є відображенням $\overline{V_p^i}$
на $U^{\,\prime}_p$ для будь-яких $i=1,2,\ldots, k, k+1.$ Отже,
$N(f, D)=N(F, D)\geqslant k+1,$ що суперечить визначенню числа $k.$
Отримана суперечність спростовує припущення, що $N(f,
\overline{D})>N(f, D).$ Лема доведена.~$\Box$
\end{proof}

\medskip
Справедливі наступні результати.

\medskip
\begin{theorem}\label{th2}
{\sl Нехай $D\subset {\Bbb R}^n,$ $n\geqslant 2,$ -- область, яка
має слабо плоску межу, а область $D^{\,\prime}\subset {\Bbb R}^n$ є
регулярною. Припустимо, $f_m,$ $m=1,2,\ldots ,$ -- послідовність
відкритих, дискретних і замкнених відображень області $D$ на
$D^{\,\prime},$ які збігаються рівномірно у $D$ до деякого
відображення $f:D\rightarrow {\Bbb R}^n$ і задовольняють
співвідношення~(\ref{eq2*A}) в кожній точці $y_0\in
\overline{D^{\,\prime}}.$

Нехай також виконана одна з двох умов:

\medskip
1) $Q\in FMO(\overline{D^{\,\prime}});$

\medskip
2) для будь-якого $y_0\in \overline{D^{\,\prime}}$ знайдеться
$\delta(y_0)>0$ таке, що для достатньо малих $\varepsilon>0$
виконуються умови
\begin{equation}\label{eq5G}
\int\limits_{\varepsilon}^{\delta(y_0)}
\frac{dt}{tq_{y_0}^{\frac{1}{n-1}}(t)}<\infty, \qquad
\int\limits_{0}^{\delta(y_0)}
\frac{dt}{tq_{y_0}^{\frac{1}{n-1}}(t)}=\infty\,.
\end{equation}
Тоді відображення $f$ є нульвимірним і замкненим у $D,$ крім того,
$f$ має неперервне продовження $\overline{f}:\overline{D}\rightarrow
\overline{D^{\,\prime}}_P,$
$\overline{f}(\overline{D})=\overline{D^{\,\prime}}_P$ і
$\overline{f}$ є нульвимірним відносно $\overline{D}_P.$}
\end{theorem}

\medskip
\begin{theorem}\label{th4}
{\sl\,Нехай за умов теореми~\ref{th2} відображення $f$ зберігає
орієнтацію. Тоді відображення $f$ є відкритим, дискретним і
замкненим у $D,$ крім того, $f$ має неперервне продовження
$\overline{f}:\overline{D}\rightarrow \overline{D^{\,\prime}}_P,$
$\overline{f}(\overline{D})=\overline{D^{\,\prime}}_P$ таке, що
$N(f, D)=N(f, \overline{D})<\infty.$ Зокрема, $\overline{f}$ є
дискретним у $\overline{D}.$ }
\end{theorem}

\medskip
Доведення теорем~\ref{th2} і~\ref{th4} випливає з лем~\ref{lem5}
і~\ref{lem6} за рахунок спеціально підібраної функції $\psi,$ див.
\cite[лема~1.3]{Sev$_3$}.~$\Box$


КОНТАКТНА ІНФОРМАЦІЯ

\medskip
\noindent{{\bf Євген Олександрович Севостьянов} \\
{\bf 1.} Житомирський державний університет ім.\ І.~Франко\\
вул. Велика Бердичівська, 40 \\
м.~Житомир, Україна, 10 008 \\
{\bf 2.} Інститут прикладної математики і механіки
НАН України, \\
вул.~Добровольського, 1 \\
м.~Слов'янськ, Україна, 84 100\\
e-mail: esevostyanov2009@gmail.com}

\medskip
\noindent{{\bf Валерій Андрійович Таргонський} \\
Житомирський державний університет ім.\ І.~Франко\\
вул. Велика Бердичівська, 40 \\
м.~Житомир, Україна, 10 008 \\
e-mail: w.targonsk@gmail.com }

\end{document}